\documentclass[11pt]{article}
\usepackage{}
\usepackage{amssymb}
\usepackage{amsfonts}
\usepackage{amsmath}
\usepackage{graphicx}
\usepackage{amsmath,amssymb,latexsym,color}
\usepackage[mathscr]{eucal}
 \makeatletter
    
    \newcommand{\Rmnum}[1]
    {\expandafter\@slowromancap\romannumeral #1@}
    \makeatother
\def\wz{\widetilde}

\newtheorem{thm}{Theorem}[section]
\newtheorem{constr}[thm]{Construction}
\newtheorem{prop}[thm]{Proposition}
\newtheorem{lemma}[thm]{Lemma}

\newtheorem{example}{Example}[section]
\newtheorem{defin}[thm]{Definition}

\newtheorem{remark}{Remark}[section]

\newcommand{\qed}{\hfill\Box\medskip}

\usepackage{CJK}
\begin{document}
\begin{CJK*}{GBK}{song}

\renewcommand{\baselinestretch}{1.3}
\title{Weakly distance-regular digraphs of valency three, \Rmnum{1}}

\author{
Yuefeng Yang,\quad Benjian Lv,\quad
Kaishun Wang\footnote{\scriptsize Corresponding author.\newline
{\em E-mail address:} yangyf@mail.bnu.edu.cn(Y.Yang), bjlv@bnu.edu.cn(B.Lv), wangks@bnu.edu.cn(K.Wang).} \\
{\footnotesize \em  Sch. Math. Sci. {\rm \&} Lab. Math. Com. Sys.,
Beijing Normal University, Beijing, 100875,  China} }
\date{}
\maketitle

\begin{abstract}
Suzuki (2004) \cite{HS04} classified thin weakly distance-regular digraphs and proposed the project to classify weakly distance-regular digraphs of
valency 3. The case of girth $2$ was classified by the third author (2004) \cite{KSW04} under the assumption of the commutativity. In this paper, we
continue this project and classify these digraphs with girth more than $2$ and two types of arcs.



\end{abstract}

\section{Introduction}

A \emph{digraph} $\Gamma$ is a pair $(X,A)$ where $X$ is a finite set of \emph{vertices} and $A\subseteq X^{2}$ is a set of \emph{arcs}. Throughout this
paper we use the term `digraph' to mean a finite directed graph with no loops. We always write $V\Gamma$ for $X$ and $A\Gamma$ for $A$. A
\emph{path} of length $r$ from $u$ to $v$ is a finite sequence of vertices $(u=w_{0},w_{1},\ldots,w_{r}=v)$
such that $(w_{t-1}, w_{t})\in A\Gamma$ for $t=1,2,\ldots,r$. A digraph is said
to be \emph{strongly connected} if, for any two distinct vertices $x$ and $y$, there is a path from $x$ to $y$. The length of a shortest path from $x$ to $y$ is called
the \emph{distance} from $x$ to $y$ in $\Gamma$, denoted by $\partial_\Gamma(x,y)$. The \emph{diameter} of $\Gamma$ is the maximum value of the distance function in $\Gamma$.
 Let $\wz{\partial}_\Gamma(x,y)=(\partial_\Gamma(x,y),\partial_\Gamma(y,x))$ and $\wz{\partial}(\Gamma)=\{\wz{\partial}_\Gamma(x,y)\mid x,y\in V\Gamma\}$.
If no confusion occurs, we  write $\partial(x,y)$ (resp. $\wz{\partial}(x,y)$) instead of $\partial_\Gamma(x,y)$ (resp. $\wz{\partial}_\Gamma(x,y)$).
An arc $(u,v)$ of $\Gamma$ is of \emph{type} $(1,r)$ if $\partial(v,u)=r$. A path $(w_{0},w_{1},\ldots,w_{r-1})$ is said to be a \emph{circuit} of length $r$ if
$\partial(w_{r-1},w_{0})=1$. A circuit is \emph{undirected} if each of its arcs is of type $(1,1)$. The \emph{girth} of $\Gamma$ is the length of a shortest circuit.

Let $\Gamma=(X,A)$ and $\Gamma'=(X',A')$ be two digraphs. $\Gamma$ and $\Gamma'$ are \emph{isomorphic} if there is a bijection $\sigma$ from $X$ to $X'$
such that $(x,y)\in A$ if and only if $(\sigma(x),\sigma(y))\in A'$. In this case, $\sigma$ is called an \emph{isomorphism} from $\Gamma$ to $\Gamma'$. An
isomorphism from $\Gamma$ to itself is called an \emph{automorphism} of $\Gamma$. The set of all automorphisms of $\Gamma$ forms a group which is called
the \emph{automorphism group} of $\Gamma$ and denoted by $\textrm{Aut}(\Gamma)$. A digraph $\Gamma$ is \emph{vertex transitive} if $\textrm{Aut}(\Gamma)$
is transitive on $V\Gamma$.

Lam \cite{CWL80} introduced a concept of distance-transitive digraphs. A connected digraph $\Gamma$ is said to be \emph{distance-transitive} if, for any
vertices $x$, $y$, $x'$ and $y'$ of $\Gamma$ satisfying $\partial(x,y)=\partial(x',y')$, there exists an automorphism $\sigma$ of $\Gamma$ such that
$x'=\sigma(x)$ and $y'=\sigma(y)$. Damerell \cite{RMD81} generalized this concept to that of distance-regular digraphs. He showed that the girth $g$ of a
distance-regular digraph of diameter $d$ is either $2$, $d$ or $d+1$, and the one with $d=g$ is a coclique extension of a distance-regular digraph with
$d=g-1$. Bannai,
Cameron and Kahn \cite{bck} proved that a distance-transitive
digraph of odd girth is a Paley tournament or a directed cycle.
Leonard and Nomura \cite{ln} proved that except directed cycles
all distance-regular digraphs with $d=g-1$ have girth $g\le 8$. In order to find `better' classes of digraphs with unbounded diameter, Damerell
\cite{RMD81} also proposed a more natural definition of distance-transitivity, i.e., weakly distance-transitivity. In \cite{KSW03}, Wang and Suzuki
introduced weakly distance-regular digraphs as a generalization of distance-regular digraphs and weakly distance-transitive digraphs.

A strongly connected digraph $\Gamma$ is said to be \emph{weakly distance-transitive}  if, for any vertices $x$, $y$, $x'$ and $y'$ satisfying
$\wz{\partial}(x,y)=\wz{\partial}(x',y')$, there exists an automorphism $\sigma$ of $\Gamma$ such that $x'=\sigma(x)$ and $y'=\sigma(y)$. A strongly
connected digraph $\Gamma$ is said to be \emph{weakly distance-regular} if, for all $\wz{h}$, $\wz{i}$, $\wz{j}\in\wz{\partial}(\Gamma)$ and
$\wz{\partial}(x,y)=\wz{h}$, the number $p_{\wz{i},\wz{j}}^{\wz{h}}:=|P_{\wz{i},\wz{j}}(x,y)|$ depends only on $\wz{h}$, $\wz{i}$, $\wz{j}$, where $P_{\wz{i},\wz{j}}(x,y)=\{z\in V\Gamma\mid\wz{\partial}(x,z)=\wz{i}~\textrm{and}~\wz{\partial}(z,y)=\wz{j}\}.$ The nonnegative integers $p_{\wz{i},\wz{j}}^{\wz{h}}$ are called the \emph{intersection numbers}. We say
that $\Gamma$ is \emph{commutative}~(resp. \emph{thin})~if $p_{\wz{i},\wz{j}}^{\wz{h}}=p_{\wz{j},\wz{i}}^{\wz{h}}$~(resp.
$p_{\wz{i},\wz{j}}^{\wz{h}}\leq1$)~for all $\wz{i}$, $\wz{j}$, $\wz{h}\in\wz{\partial}(\Gamma)$. Note that a weakly distance-transitive digraph is weakly
distance-regular.

Let $G$ be a finite group and $S$ a subset of $G$ not containing the identity. The \emph{Cayley digraph} $\Gamma=\textrm{Cay}(G,S)$ is a digraph with the vertex set $G$ and the arc set $\{(x,sx)\mid x\in G,~s\in S\}.$

In \cite{KSW03}, Wang and Suzuki determined all commutative $2$-valent weakly distance-regular digraphs. In \cite{HS04}, Suzuki determined all thin weakly
distance-regular digraphs and proved the nonexistence of noncommutative weakly distance-regular digraphs of valency $2$. Moreover, he proposed the project
to classify weakly distance-regular digraphs of valency $3$. In \cite{KSW04}, Wang classified all commutative weakly distance-regular digraphs of valency
$3$ and girth $2$. In this paper, we continue this project, and obtain the following result.

\begin{thm}\label{Main}
Let $\Gamma$ be a weakly distance-regular digraph of valency $3$ and girth more than $2$. If $\Gamma$ has two types of arcs, then $\Gamma$ is
isomorphic to one of the following digraphs:
\begin{itemize}
\item [{\rm(i)}] ${\rm Cay}(\mathbb{Z}_{4}\times\mathbb{Z}_{g},\{(0,1),(2,1),(1,0)\})$, where $g=3$ or $g\geq5$.\vspace{-0.3cm}

\item [{\rm(ii)}] $\Gamma_{q,2mq,1}$, $\Gamma_{q,mq+2,q}$ or $\Gamma_{q,2mq-2q+2t,q+1-t}$ in Construction \ref{construcution}, where $q\geq 3$, $m\geq1$ and $2\leq t\leq q-1$.
\end{itemize}
\end{thm}

This paper is organized as follows. In Section 2, we construct two families of weakly distance-regular digraphs of valency $3$. In Section 3, we discuss
some properties for circuits of weakly distance-regular digraphs. In Section 4, we prove our main theorem.

\section{Constructions}

In this section, we construct two families of weakly distance-regular digraphs of valency $3$. For any element $x$ in a residue class ring, we always assume that $\hat{x}$ denotes the minimum nonnegative integer in $x$. Denote
$\beta(w)=(1+(-1)^{w+1})/2$ for any integer $w$.

\begin{prop}\label{G_Wdrdg}
Let $g\geq3$. Then $\Gamma_{g}:={\rm Cay}(\mathbb{Z}_{4}\times \mathbb{Z}_{g},\{(1,0),(0,1),(2,1)\})$ is a weakly distance-regular digraph if and only if $g\neq4$.
\end{prop}
\textbf{Proof.}    For any vertex $(a,b)$ distinct with $(0,0)$,
we have
\begin{eqnarray}
\wz{\partial}((0,0),(a,b))=\left\{
\begin{array}{ll}
(\hat{a},4-\hat{a}), & \textrm{if}\ b=0,\\
(\hat{b}+\beta(\hat{a}),g-\hat{b}+\beta(\hat{a})), & \textrm{if}\ b\neq0.
\end{array} \right.\nonumber
\end{eqnarray}

Suppose $g\neq4$. We will show that $\Gamma_{g}$ is weakly distance-transitive. Let $(a,b)$ and $(x,y)$ be any two vertices satisfying
$\wz{\partial}((0,0),(a,b))=\wz{\partial}((0,0),(x,y))$. It suffices to verify that there exists an automorphism $\sigma$ of $\Gamma_{g}$ such that
$\sigma(0,0)=(0,0)$ and $\sigma(a,b)=(x,y)$. If $(a,b)=(x,y)$, then the identity permutation is a desired automorphism. Now suppose $(a,b)\neq(x,y)$. Then $b\neq0$, $y\neq0$ and
$
(\hat{b}+\beta(\hat{a}),g-\hat{b}+\beta(\hat{a}))=(\hat{y}+\beta(\hat{x}),g-\hat{y}+\beta(\hat{x})).
$
It follows that $b=y$ and $a-x=2$. Let $\sigma$ be the permutation on $V\Gamma_{g}$ such that
\begin{eqnarray}
\sigma(x,y)=\left\{
\begin{array}{ll}
(x,y), & \textrm{if}\ y\neq b,\\
(x+2,y), & \textrm{if}\ y=b.
\end{array} \right.\nonumber
\end{eqnarray}
Routinely,  $\sigma$ is a desired automorphism.

In $\Gamma_{4}$,   $\wz{\partial}((0,0),(0,2))=\wz{\partial}((0,0),(2,0))=(2,2)$. But
$P_{(1,3),(3,3)}((0,0),(0,2))=\{(1,0)\}$ and $P_{(1,3),(3,3)}((0,0),(2,0))=\emptyset.$
Hence, $\Gamma_{4}$ is not a weakly distance-regular digraph.
$\qed$

\begin{constr}\label{construcution}
Let $q$, $s$, $k$ be integers with $q>2$,~$s>2$ and $\max\{1,q-s+2\}\leq k\leq q$. Write $s=2mq+p$ with $m\geq0$ and $0\leq p<2q$. Let $\Gamma_{q,s,k}$
be the digraph with the vertex set $\mathbb{Z}_{q}\times\mathbb{Z}_{s}$ whose arc set consists of $((a,b),(a+1,b))$, $((a,c),(a,c+1))$,
$((a,d),(a+1,d-1))$, $((a,s-1),(a-k+1,0))$ and $((a,0),(a+k,s-1))$, where $c\neq s-1$ and $d\neq0$. See Figure 1.
\end{constr}

\begin{figure}[hbt]
  \centering
   \scalebox{0.8}{\includegraphics{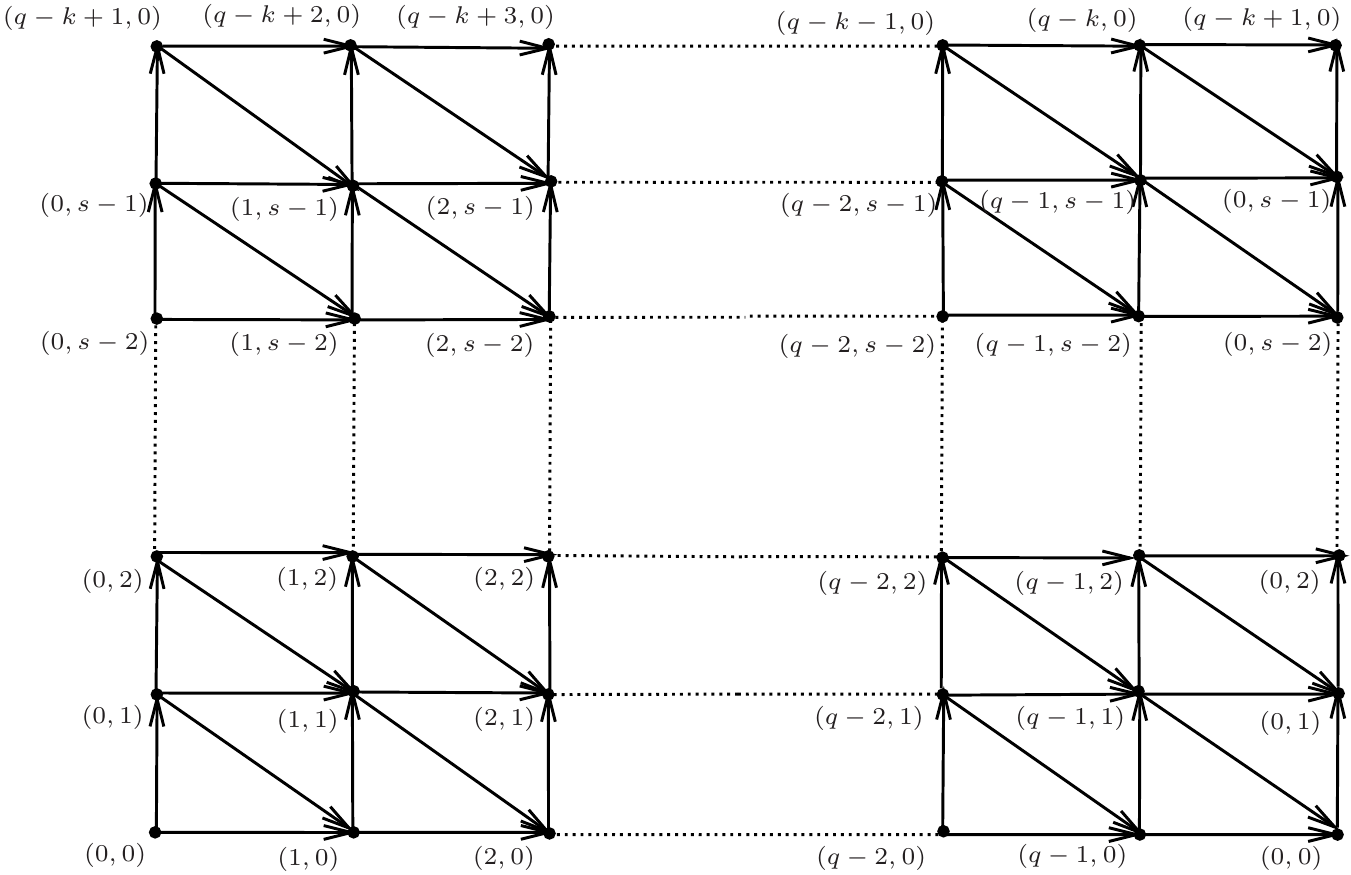}}
   \caption{The digraph $\Gamma_{q,s,k}$.}
\end{figure}

In the following, we will prove that $\Gamma_{q,s,k}$ is a weakly distance-regular digraph if and only if one of the following holds:

C1: $p=0$ and $k=1$.

C2: $p=q+2$ or $p=2$, and $k=q$.

C3: $4\leq p\leq 2q-2$, $p$ is even and $k=q+1-p/2$.

\begin{lemma}\label{QSK_VT}
$\Gamma_{q,s,k}$ is a vertex transitive digraph.
\end{lemma}
\textbf{Proof.}~Pick any vertex $(a,b)$. It  suffices to show that there exists an automorphism $\sigma$ of $\Gamma_{q,s,k}$ such that
$\sigma(0,0)=(a,b)$. Let $\sigma$ be the permutation on $V\Gamma_{q,s,k}$ such that
\begin{eqnarray}
\sigma(x,y)=\left\{
\begin{array}{ll}
(x+a,y+b), & \textrm{if}\ \hat{y}\in\{0,1,2,\ldots,s-1-\hat{b}\},\\
(x+a-k+1,y+b), & \textrm{otherwise.}
\end{array} \right.\nonumber
\end{eqnarray}
Routinely,   $\sigma$ is a desired automorphism. $\qed$

For any two integers $i$ and $j$, we always write $i\equiv j$ instead of $i\equiv j~(\textrm{mod}~q)$. For any vertex $(a,b)$ of $\Gamma_{q,s,k}$, let $f(a,b)$, $g(a,b)$ and
$h(a)$ be nonnegative integers less than $q$ such that
\begin{eqnarray}
f(a,b)\equiv \hat{a}+\hat{b}-k-p+1,~g(a,b)\equiv q-\hat{a}-\hat{b}~\textrm{and}~h(a)\equiv k-\hat{a}-1.\nonumber
\end{eqnarray}
By the structure of $\Gamma_{q,s,k}$, we have
\begin{eqnarray}
\wz{\partial}((0,0),(a,b))=(\min\{\hat{a}+\hat{b},s-\hat{b}+f(a,b)\},\min\{\hat{b}+g(a,b),s-\hat{b}+h(a)\}).\nonumber
\end{eqnarray}

\begin{lemma}\label{QSK_Dis}
Let {\rm C1, C2} or {\rm C3} hold. In $\Gamma_{q,s,k}$, $\partial((0,0),(a,b))=\hat{a}+\hat{b}$ if and only if $\partial((a,b),(0,0))=\hat{b}+g(a,b)$.
\end{lemma}
\textbf{Proof.} Let $M=s-2\hat{b}-\hat{a}+f(a,b)$ and $N=s-2\hat{b}+h(a)-g(a,b)$. We only need to prove
$M>0$ if and only if $N>0.$
Note that $f(a,b)+g(a,b)$ equals to $k-1$ or $q+k-1$ and $h(a)$ equals to $k-\hat{a}-1$ or $q+k-\hat{a}-1$.

\textbf{Case 1}.~$f(a,b)+g(a,b)=k-1$ and $h(a)=k-\hat{a}-1$, or $f(a,b)+g(a,b)=q+k-1$ and $h(a)=q+k-\hat{a}-1$.

In this case, it is routine to check $M=N$, as desired.

\textbf{Case 2}.~$f(a,b)+g(a,b)=k-1$ and $h(a)=q+k-\hat{a}-1$.

Note that $M=N-q$ and $k\neq q$. Therefore,  there exists an even number $n$ such that $s+2k-2=nq$. If $M>0$, then $N>0$. Conversely, suppose $N>0$. Let $\hat{b}=n'q+r'$
with $0\leq r'<q$. If $g(a,b)=q-\hat{a}-r'$, then $k-1+r'<\hat{a}+r'<q$, which implies that
$N=s+2k-2-2n'q-k+1-r'>q$. If $g(a,b)=2q-\hat{a}-r'$, then $f(a,b)=k-1+\hat{a}+r'-2q$. Hence, $q<k-1+r'<2q$, which implies that $N=s+2k-2-2n'q-2q+(q+1)-k-r'>q$.
Thus, $M>0$ and the desired result holds.

\textbf{Case 3}.~$f(a,b)+g(a,b)=q+k-1$ and $h(a)=k-\hat{a}-1$.

Similar to Case 2, the desired result follows.$\qed$

By Lemma \ref{QSK_VT}, for  vertices $(a,b)$ and $(x,y)$ of $\Gamma_{q,s,k}$, we have
\begin{eqnarray}
\wz{\partial}((a,b),(x,y))=\left\{
\begin{array}{ll}
\wz{\partial}((0,0),(x-a,y-b)), & \textrm{if}\ \hat{y}\in\{\hat{b},\hat{b}+1,\ldots,s-1\},\\
\wz{\partial}((0,0),(x-a+k-1,y-b)), & \textrm{otherwise.}
\end{array} \right.\nonumber
\end{eqnarray}

\begin{prop}\label{QSK_Wdrdgiff}
$\Gamma_{q,s,k}$ is a weakly distance-regular digraph if and only if {\rm C1, C2} or {\rm C3} holds.
\end{prop}
\textbf{Proof.}~``$\Longleftarrow$" We will prove that $\Gamma_{q,s,k}$ is weakly distance-transitive. Let $(a,b)$ and $(x,y)$ be two vertices satisfying
$\wz{\partial}((0,0),(a,b))=\wz{\partial}((0,0),(x,y))$. It suffices to find $\sigma\in\textrm{Aut}(\Gamma_{q,s,k})$ such that $\sigma(0,0)=(0,0)$ and
$\sigma(a,b)=(x,y)$.

\textbf{Case 1}.~$\partial((0,0),(a,b))=\hat{a}+\hat{b}$.

Suppose $\partial((0,0),(x,y))=\hat{x}+\hat{y}$. Then $g(a,b)=g(x,y)$. By Lemma \ref{QSK_Dis}, we have $\hat{b}+g(a,b)=\hat{y}+g(x,y)$. This implies that
$a=x$ and $b=y$. Hence, the identity permutation is a desired automorphism.

Suppose $\partial((0,0),(x,y))=s-\hat{y}+f(x,y)$. Then $\hat{x}\equiv\hat{a}+\hat{b}+k-1\equiv f(a,b)$. Hence, $\hat{x}=f(a,b)$ and $g(a,b)=h(x)$. By
Lemma \ref{QSK_Dis}, we have $\hat{b}+g(a,b)=s-\hat{y}+h(x)$. This implies $\hat{y}=s-\hat{b}$. Let $\sigma$ be the permutation on $V\Gamma_{q,s,k}$ such
that
\begin{eqnarray}
\sigma(a,b)=\left\{
\begin{array}{ll}
(a,b), & \textrm{if}\ b=0,\\
(f(a,b),-b), & \textrm{if}\ b\neq0.
\end{array} \right.\nonumber
\end{eqnarray}
Routinely, $\sigma$ is a desired automorphism.

\textbf{Case 2}.~$\partial((0,0),(a,b))=s-\hat{b}+f(a,b)$.

Suppose $\partial((0,0),(x,y))=s-\hat{y}+f(x,y)$. Then $\hat{y}-\hat{b}=f(x,y)-f(a,b)$. We have $\hat{y}-\hat{b}\equiv\hat{x}+\hat{y}-\hat{a}-\hat{b}$.
This implies $x=a$. By Lemma \ref{QSK_Dis}, one gets $s-\hat{b}+h(a)=s-\hat{y}+h(x)$, which implies that $y=b$. Hence, the identity permutation is a desired
automorphism.

Suppose $\partial((0,0),(x,y))=\hat{x}+\hat{y}$. It is similar to Case 1 and the desired result holds.

``$\Longrightarrow$"~Suppose C1, C2 and C3 do not hold. Let $e=(0,0)$, $z=(0,1)$, $w=(k,s-1)$, $t=\lfloor p/q \rfloor$ and $\alpha(v)=(3+(-1)^{v})/4$ for
$v\in\mathbb{Z}$.

\textbf{Case 1}.~$k\neq q$ and $2((-1)^{t}\alpha(p)+1)+qt<2k+p\leq 2(q-\alpha(p)+1)$.

Let $x=(0,\alpha(s)+s/2)$ and $y=(q+\alpha(p)+1-k-p/2,(m-1)q+p+k-1)$. In this case, $\wz{\partial}(e,x)=\wz{\partial}(e,y)$. But $z\in
P_{(1,q),\wz{\partial}(z,x)}(e,x)$ and $P_{(1,q),\wz{\partial}(z,x)}(e,y)=\emptyset.$

\textbf{Case 2}.~$k\neq q$, and $2k+p\leq 2((-1)^{t}\alpha(p)+1)+qt$ or $2(q-\alpha(p)+2)\leq 2k+p$.

Let $x=(k,\alpha(s)-1+s/2)$ and $y=(k,\alpha(s)+s/2)$. In this case, $\wz{\partial}(e,x)=\wz{\partial}(e,y)$. But $w\in
P_{(1,q),\wz{\partial}(w,x)}(e,x)$ and $P_{(1,q),\wz{\partial}(w,x)}(e,y)=\emptyset.$

\textbf{Case 3}.~$k=q$ and $3\leq p\leq q+1$.

Let $x=(q-2,mq+2)$ and $y=(q-2,mq+1)$. In this case, $\wz{\partial}(e,x)=\wz{\partial}(e,y)$. But $z\in P_{(1,q),\wz{\partial}(z,x)}(e,x)$ and $P_{(1,q),\wz{\partial}(z,x)}(e,y)=\emptyset.$

\textbf{Case 4}.~$k=q$, and $p\leq1$ or $q+3\leq p$.

Let $x=(q-1,mq-tq+p)$ and $y=(0,mq+tq)$. In this case, $\wz{\partial}(e,x)=\wz{\partial}(e,y)$. But $z\in P_{(1,q),\wz{\partial}(z,x)}(e,x)$
and $P_{(1,q),\wz{\partial}(z,x)}(e,y)=\emptyset.$

In all above cases, $\Gamma_{q,s,k}$ is not weakly distance-regular and the desired result holds.$\qed$

Finally, we shall show that every weakly distance-regular digraph $\Gamma_{q,s,k}$ is a Cayley digraph.

\begin{prop}
Let $d=\frac{p}{2(q,p)}$, $l=\max\{w\mid2^{w}~\textrm{divides}~(q,p)\}$, $h=\frac{s}{2^{l}}$, $i=2\{d\}$ and $u$ be an integer such that $2^{i}q$ divides
$(up-(q,p))$, where $\{d\}$ denotes the fractional part of $d$ and $(q,p)$ denotes the greatest common divisor of $q$ and $p$. Then the weakly
distance-regular digraph $\Gamma_{q,s,k}$ is isomorphic to one of the following Cayley digraphs:
\begin{itemize}
\item[${\rm(i)}$] ${\rm Cay}(\mathbb{Z}_{q}\times\mathbb{Z}_{2mq},\{(1,0),(0,1),(1,2mq-1)\})$, $m\geq1$ and $q\geq3$.\vspace{-0.3cm}

\item[${\rm(ii)}$] ${\rm Cay}(\mathbb{Z}_{(mq+2)q},\{1,mq+2,mq+1\})$, $m\geq1$ and $q\geq3$.\vspace{-0.3cm}

\item[${\rm(iii)}$] ${\rm Cay}(\mathbb{Z}_{2^{i}q}\times\mathbb{Z}_{2^{-i}(2mq+p)},\{(2^{i},ih),(2^{i}ud,1),(2^{i}-2^{i}ud,ih-1)\})$, where $q\geq 3$,
    $m\geq0$, $4\leq p\leq 2q-2$ and $p$ is even.
\end{itemize}
\end{prop}
\textbf{Proof.}~If C1 holds, then (i) is obvious. If C2 holds, then the mapping $\sigma$ from $\Gamma_{q,s,k}$ to the digraph in (ii) satisfying
$\sigma(a,b)=\hat{a}(mq+2)+\hat{b}$ is an isomorphism.

Now suppose C3 holds. Let $\sigma$ be the mapping from $\Gamma_{q,s,k}$ to the digraph in (iii) such that
$\sigma(a,b)=(2^{i}\hat{a}+2^{i}ud\hat{b},ih\hat{a}+\hat{b})$. Note that $\sigma$ is well defined. We will show that $\sigma$ is injective. It is clear
for $i=0$. If $i=1$, then $l\geq 1.$ Assume that $\sigma(x_{1},y_{1})=\sigma(x_{2},y_{2})$ for $(x_{1},y_{1}),(x_{2},y_{2})\in V\Gamma_{q,s,k}$. Let
$x=2ud(\widehat{y_{2}}-\widehat{y_{1}})-2(\widehat{x_{1}}-\widehat{x_{2}})$ and $y=(\widehat{y_{2}}-\widehat{y_{1}})-h(\widehat{x_{1}}-\widehat{x_{2}})$.
We have $2q| x$ and $(mq+p/2)|y$.
Hence, $h| (\widehat{y_{2}}-\widehat{y_{1}})$. We claim $2^{j}|(\widehat{y_{2}}-\widehat{y_{1}})$ for $1\leq j\leq l$. Note that
$2|(\widehat{y_{2}}-\widehat{y_{1}})$. Suppose $2^{j}|(\widehat{y_{2}}-\widehat{y_{1}})$ for some $j<l$. Since $2^{j}| y$, one gets $2^{j}|
(\widehat{x_{1}}-\widehat{x_{2}})$ and $2^{j+1}| x$, which imply that $2^{j+1}|(\widehat{y_{2}}-\widehat{y_{1}})$. So our claim is valid. By $(2^{l},h)=1$, we obtain $(2mq+p)| (\widehat{y_{2}}-\widehat{y_{1}})$. Thus, $y_{1}=y_{2}$ and $x_{1}=x_{2}$. Therefore $\sigma$ is a bijection. One can verify that
$((x_{1},y_{1}),(x_{2},y_{2}))$ is an arc if and only if $(\sigma(x_{1},y_{1}),\sigma(x_{2},y_{2}))$ is an arc. Hence, $\sigma$ is an isomorphism.$\qed$

\section{Circuits}

In this section, we will discuss some properties for circuits of weakly distance-regular digraphs.

Let $\Gamma$ be a digraph. Let $R=\{\Gamma_{\wz{i}}\mid\wz{i}\in\wz{\partial}(\Gamma)\}$, where $\Gamma_{\wz{i}}=\{(x,y)\in V\Gamma\times V\Gamma\mid\wz{\partial}(x,y)=\wz{i}\}$. If $\Gamma$ is weakly distance-regular, then $(V\Gamma,R)$ is an association scheme. For more information about association schemes,
see~\cite{EB84,PHZ96}. For two nonempty subsets $E$, $F\subseteq R$, define
\begin{eqnarray}
EF:=\{\Gamma_{\wz{h}}\mid\sum_{\Gamma_{\wz{i}}\in E}\sum_{\Gamma_{\wz{j}}\in F}p_{\wz{i},\wz{j}}^{\wz{h}}\neq0\}, \nonumber
\end{eqnarray}
and write $\Gamma_{\wz{i}}\Gamma_{\wz{j}}$ instead of $\{\Gamma_{\wz{i}}\}\{\Gamma_{\wz{j}}\}$. For each nonempty subset $F$ of $R$, define $\langle
F\rangle$ to be the minimal equivalence relation containing $F$. Let
\begin{eqnarray}
V\Gamma/F:=\{F(x)\mid x\in V\Gamma\}\quad\textrm{and}\quad \Gamma_{\wz{i}}^{F}:=\{(F(x),F(y))\mid y\in F\Gamma_{\wz{i}}F(x)\},\nonumber
\end{eqnarray}
where $F(x):=\{y\in V\Gamma\mid(x,y)\in \cup_{f\in F}f\}$. The digraph $(V\Gamma/F,\ \cup_{(1,s)\in\wz{\partial}(\Gamma)}\Gamma_{1,s}^{F})$ is said to be
the \emph{quotient digraph} of $\Gamma$ over $F$, denoted by $\Gamma/F$. The size of $\Gamma_{\wz{i}}(x):=\{y\in V\Gamma\mid\wz{\partial}(x,y)=\wz{i}\}$
depends only on $\wz{i}$, denoted by $k_{\wz{i}}$. For any $(a,b)\in\wz{\partial}(\Gamma)$, we usually write $k_{a,b}$ (resp. $\Gamma_{a,b}$) instead of $k_{(a,b)}$ (resp. $\Gamma_{(a,b)}$).

Now we shall introduce some basic results which are used frequently in this paper.
\begin{lemma}\label{yinyong}
Let $\Gamma$ be a weakly distance-regular digraph. For each $\wz{i}:=(a,b)\in\wz{\partial}(\Gamma)$, define $\wz{i}^{*}=(b,a)$.\vspace{-0.3cm}
\begin{itemize}
\item[${\rm(i)}$] $k_{\wz{h}}p_{\wz{i},\wz{j}}^{\wz{h}}=k_{\wz{i}}p_{\wz{h},\wz{j}^{*}}^{\wz{i}}=k_{\wz{j}}p_{\wz{i}^{*},\wz{h}}^{\wz{j}}$.\vspace{-0.3cm}

\item[${\rm(ii)}$] $k_{\wz{i}}k_{\wz{j}}=\sum_{\wz{h}\in\wz{\partial}(\Gamma)}k_{\wz{h}}p_{\wz{i},\wz{j}}^{\wz{h}}$.\vspace{-0.3cm}

\item[${\rm(iii)}$] $|\Gamma_{\wz{i}}\Gamma_{\wz{j}}|\leq(k_{\wz{i}},k_{\wz{j}})$.
\end{itemize}
\end{lemma}
\textbf{Proof.}~See Proposition 2.2 in \cite[pp. 55-56]{EB84} and \cite[Proposition 5.1]{ZA99}.$\qed$

In the remaining of this paper, we always assume that $\Gamma$ is a weakly distance-regular digraph of valency $3$ satisfying $k_{1,q-1}=1$ and
$k_{1,g-1}=2$, where $q,g\geq3$ and $q\neq g$. Let $A_{i,j}$ denote a binary matrix with rows and columns indexed by $V\Gamma$ such that $(A_{i,j})_{x,y}=1$ if and only if $\wz{\partial}(x,y)=(i,j).$

\begin{lemma}\label{Arc} The following hold:
\begin{eqnarray}
A_{1,q-1}A_{1,g-1}=A_{1,g-1}A_{1,q-1},\label{arci}\\
A_{1,g-1}A_{g-1,1}=A_{g-1,1}A_{1,g-1}.\label{arcii}
\end{eqnarray}
\end{lemma}
\textbf{Proof.}~By Lemma \ref{yinyong} (iii), we may assume that
$$
A_{1,g-1}A_{1,q-1}=A_{i,j}\quad{\rm and}\quad A_{1,q-1}A_{1,g-1}=A_{s,t},\quad i,s\in\{1,2\}.
$$
 We claim that $i=s=2$. Suppose $i=1$. Then $j=g-1$ because of $k_{1,q-1}=1$. By Lemma \ref{yinyong} (i), we get $p_{(g-1,1),(1,g-1)}^{(1,q-1)}=2p_{(1,g-1),(1,q-1)}^{(1,g-1)}=2$.
By Lemma \ref{yinyong} (iii), $A_{g-1,1}A_{1,g-1}=2I+2A_{1,q-1}$, contrary to the fact that $A_{g-1,1}A_{1,g-1}$ is a symmetric matrix. Hence, $i=2$. Similarly, $s=2$ and our claim is valid.

Pick a path $(x_{0},x_{1},x_{2})$ with $\wz{\partial}(x_{0},x_{1})=(1,g-1)$ and $\wz{\partial}(x_{1},x_{2})=(1,q-1)$. Then $\partial(x_{2},x_{0})=j$. We may choose a path $(x_{2}, x_{3},\ldots,x_{j+1},x_{0})$. Since $\Gamma$ has just two types of arcs, there exists an $i\in\{1,2,\ldots,j+1\}$ such that $\wz{\partial}(x_{i},x_{i+1})=(1,q-1)$ and $\wz{\partial}(x_{i+1},x_{i+2})=(1,g-1)$, where $x_{j+2}=x_{0}$ and $x_{j+3}=x_{1}$. Since $\wz{\partial}(x_{i},x_{i+2})=(2,t)$, one has $t\leq j$. Similarly, $j\leq t$. Hence, $j=t$ and the (\ref{arci}) holds.

In view of Lemma \ref{yinyong} (iii), we have
\begin{eqnarray}
A_{1,g-1}A_{g-1,1} & = & 2I+p_{(1,g-1),(g-1,1)}^{(s,s)}~A_{s,s},\quad s\geq2,\label{Arc(1)}\\
A_{g-1,1}A_{1,g-1} & = & 2I+p_{(g-1,1),(1,g-1)}^{(t,t)}~A_{t,t},\quad t\geq2.\label{Arc(2)}
\end{eqnarray}
By Lemma \ref{yinyong} (ii), we have $k_{s,s}p_{(1,g-1),(g-1,1)}^{(s,s)}=k_{t,t}p_{(g-1,1),(1,g-1)}^{(t,t)}=2$, which implies that $p_{(1,g-1),(g-1,1)}^{(s,s)},p_{(g-1,1),(1,g-1)}^{(t,t)}\in\{1,2\}$. Let $x_{0}$ and $x_{s}$ be two vertices satisfying $\wz{\partial}(x_{0},x_{s})=(s,s)$. Suppose $p_{(1,g-1),(g-1,1)}^{(s,s)}=2$.  Pick two distinct vertices $x,y\in P_{(1,g-1),(g-1,1)}(x_{0},x_{s})$. By (\ref{Arc(2)}), $\wz{\partial}(x,y)=(t,t)$. It follows that $p_{(g-1,1),(1,g-1)}^{(t,t)}=2$. Similarly, if $p_{(g-1,1),(1,g-1)}^{(t,t)}=2$, then $p_{(1,g-1),(g-1,1)}^{(s,s)}=2$ by (\ref{Arc(1)}). Hence, $p_{(1,g-1),(g-1,1)}^{(s,s)}=p_{(g-1,1),(1,g-1)}^{(t,t)}$.
In order to show (\ref{arcii}), we shall  prove $s=t$. Pick $x\in
P_{(1,g-1),(g-1,1)}(x_{0},x_{s})$ and a path $P:=(x_{0},\ x_{1},\ldots,x_{s})$.

\textbf{Case 1}. $P$ contains an arc of type $(1,g-1)$.

By (\ref{arci}), without loss of generality, we may assume that $\wz{\partial}(x_{0},x_{1})=(1,g-1)$. Pick $y\in\Gamma_{1,g-1}(x_{s})\setminus\{x\}$. In view of (\ref{Arc(2)}), if $x\neq x_{1}$, then $\partial(x_{1},x)=t\leq s$; if $x=x_{1}$, then $\partial(x,y)=t\leq s$.

\textbf{Case 2}. $P$ only contains arcs of type $(1,q-1)$.

In this case, $A_{1,q-1}^{s}\neq I$. By (\ref{arci}), there exists a path $(x_{0},y_{1},y_2,\ldots,y_{s},x)$ containing the unique arc $(x_{0},y_{1})$ of type $(1,g-1)$. If $x=y_{1}$, by Lemma \ref{yinyong} (iii),  we have $A_{1,q-1}^{s}=I$, a contradiction. Therefore, $x\neq y_{1}$. By (\ref{Arc(2)}), one has $\partial(y_{1},x)=t\leq s$.

Similarly, $t\geq s$, which implies $s=t$, as desired.$\qed$

In the following, let $F=\langle \Gamma_{1,g-1}\rangle$ and fix $x\in V\Gamma$. Then $\Gamma/F$ is isomorphic to a circuit $C_{m}$ of length $m$. Let $\Delta$ be a digraph with the vertex set $F(x)$ such that $(y,z)$ is an arc of $\Delta$ if $(y,z)$ is an arc of type $(1,g-1)$ in $\Gamma$.

\begin{lemma}\label{f(x)}
Suppose that every circuit of length $g$ contains arcs of the same type in $\Gamma$. Then $\Delta_{t,g-t}(x_{0})=\Gamma_{t,g-t}(x_{0})$ for each $x_{0}\in F(x)$ and $t\in\{1,2,\ldots,g-1\}$.
\end{lemma}
\textbf{Proof.}~Note that every arc of type $(1,g-1)$ is contained in a circuit of length $g$ with all arcs of type $(1,g-1)$. It follows that,
for any such circuit $(x_0,x_1,\ldots, x_{g-1})$, we have  $\wz{\partial}_{\Gamma}(x_0, x_i)=(i, g-i)$, where $1\leq i\leq g-1.$
Then every arc of $\Delta$ is contained in a circuit of length $g$ in $\Delta$.

For any $x_{t}\in \Gamma_{t,g-t}(x_{0})$, there exists a circuit $C_{g}:=(x_{0},x_{1},\ldots,x_{t},\ldots,x_{g-1})$ in $\Gamma$. Hence, $C_{g}$ only contains the arcs of same type. Suppose that each arc of $C_{g}$ is of type $(1,q-1)$. Then, $q<g$ and every circuit of length $q$ in $\Gamma$ only contains arcs of type $(1,q-1)$. It follows that $A_{1,q-1}^{q}=I$. Since $x_{0}\neq x_{l}$ for $1\leq l\leq g-1$, $k_{1,q-1}=1$ implies that $g$ is the minimum positive integer such that $A_{1,q-1}^{g}=I$, a contradiction. Consequently, each arc of $C_{g}$ is of type $(1,g-1)$. Therefore, $(x_{0},x_{t})\in\Delta_{t,g-t}$; and so $\Gamma_{t,g-t}(x_{0})\subseteq\Delta_{t,g-t}(x_{0})$. Conversely, pick any $x_{t}\in \Delta_{t,g-t}(x_{0})$. Then, in $\Gamma$, there exists a circuit $(x_{0},x_{1},\ldots,x_{t},\ldots,x_{g-1})$ each of whose arcs is of type $(1,g-1)$. Hence, $(x_{0},x_{t})\in \Gamma_{t,g-t}$; and so $\Delta_{t,g-t}(x_{0})\subseteq\Gamma_{t,g-t}(x_{0})$. Thus, the desired result holds. $\qed$

\begin{lemma}\label{F(x)VGammaArcGG}
If $F(x)=V\Gamma$, then there exists a circuit of length $g$ containing different types of arcs.
\end{lemma}
\textbf{Proof.}~Suppose for the contrary that every circuit of length $g$ contains the same type of arcs. By the Lemma \ref{f(x)}, $\Gamma_{t,g-t}=\Delta_{t,g-t}$ for any $1\leq t\leq g-1$. By (\ref{arcii}), the proof of Proposition 4.3 in \cite{KSW03} implies that $\Delta$ is isomorphic to
$\Gamma_{1}:=\textrm{Cay}(\mathbb{Z}_{2g},\{1,g+1\})$ or $ \Gamma_{2}:=\textrm{Cay}(\mathbb{Z}_{g}\times\mathbb{Z}_{g},\{(0,1),(1,0)\}).$

\textbf{Case 1}.~$\Delta\simeq\Gamma_{1}$.

Choose $y\in\mathbb{Z}_{2g}\setminus\{0,1,g+1\}$ and $t\in\mathbb{Z}_{2g}$ such that $\wz{\partial}_{\Gamma}(0,y)=(1,q-1)$, $\hat{t}\equiv \hat{y}~(\textrm{mod}~g)$ and $\hat{t}\in\{0,2,3,\ldots,g-1\}$. Since $(y+1,y+2,\ldots,y-t+g-1,0,y)$ is a path of length $g-\hat{t}$, $\partial_{\Gamma}(y+1,y)=g-1\leq g-\hat{t}$. It follows that $t=0$, and so $y=g$. Therefore, $\wz{\partial}_{\Gamma}(0,g)=(1,q-1)$. Similarly, $\wz{\partial}_{\Gamma}(g,0)=(1,q-1)$. Hence, $q=2$, a contradiction.

\textbf{Case 2}.~$\Delta\simeq\Gamma_{2}$.

Pick $(i,j)\in\Gamma_{1,q-1}(0,0)$. Since $\wz{\partial}_{\Delta}((0,0)(0,j))=(\hat{j},g-\hat{j})$, by Lemma \ref{f(x)}, we have $\wz{\partial}_{\Gamma}((0,0)(0,j))=(\hat{j},g-\hat{j})$. It follows that $i\neq 0$. By Lemma \ref{yinyong} (i), one gets
$
p_{(\hat{i},g-\hat{i}),(\hat{j},g-\hat{j})}^{(1,q-1)}=
k_{\hat{i},g-\hat{i}}~p_{(1,q-1),(g-\hat{j},\hat{j})}^{(\hat{i},g-\hat{i})}.$
Since $(i,j)\in P_{(1, q-1),(g-\hat{j}, \hat{j})}((0,0),(i,0))$ in $\Gamma$,  $p_{(1,q-1),(g-\hat{j},\hat{j})}^{(\hat{i},g-\hat{i})}=1$, which implies that $p_{(\hat{i},g-\hat{i}),(\hat{j},g-\hat{j})}^{(1,q-1)}=k_{\hat{i},g-\hat{i}}$.

Let $((a,b),(a',b'))$ be
an arc of type $(1,q-1)$. Then
$
P_{(\hat{i},g-\hat{i}),(\hat{j},g-\hat{j})}((a,b),(a',b'))=\Gamma_{\hat{i},g-\hat{i}}(a,b).$
Since $(a+i,b),(a,b+i)\in\Delta_{\hat{i},g-\hat{i}}(a,b)$, by Lemma \ref{f(x)}, $(a',b')\in\Gamma_{\hat{j},g-\hat{j}}(a+i,b)\cap\Gamma_{\hat{j},g-\hat{j}}(a,b+i)$. By Lemma \ref{f(x)} again, $(a',b')\in\{(a+i+j,b),(a+i,b+j)\}\cap\{(a+j,b+i),(a,b+i+j)\}$. Since $i\neq0$, we have $(a',b')=(a+i,b+j)=(a+j,b+i)$, which implies that $i=j$.
Thus,
$\Gamma\simeq\textrm{Cay}(\mathbb{Z}_{g}\times\mathbb{Z}_{g},\{(1,0),(0,1),(i,i)\})$. Since $g\neq q$, one gets $i\neq1$ and
$\wz{\partial}_{\Gamma}((0,0),(1,1))=\wz{\partial}_{\Gamma}((0,0),(i,i+1))$. But $(1,0)\in P_{(1,g-1),(1,g-1)}((0,0),(1,1)) $ and $
P_{(1,g-1),(1,g-1)}((0,0),(i,i+1))=\emptyset$ in $\Gamma$,
a contradiction.
$\qed$

\begin{lemma}\label{Circuit}
Every circuit of length $q$ in $\Gamma$ only contains the arcs of the same type. In particular,
\begin{eqnarray}\label{2,g-2}
A_{1,q-1}^{2}=A_{2,q-2}.
\end{eqnarray}
\end{lemma}
\textbf{Proof.} If $F(x)=V\Gamma$, then $q<g$ by Lemma \ref{F(x)VGammaArcGG} and the desired result follows. Suppose $F(x)\neq V\Gamma$. Assume the contrary, namely, there exists a circuit $(x_{0},x_{1},\ldots,x_{q-1})$ containing arcs of different types. Since $\Gamma/F\simeq C_{m}$ with $m\geq2$, there exist at least two arcs of type $(1,q-1)$ in this circuit. By (\ref{arci}), we may assume that $\wz{\partial}(x_{0},x_{1})=\wz{\partial}(x_{1},x_{2})=(1,q-1)$ and $\wz{\partial}(x_{q-1},x_{0})=(1,g-1)$. By the claim in Lemma~\ref{Arc}, $\wz{\partial}(x_{q-1},x_{1})=(2,q-2)$. Since $k_{1,q-1}=1$, by Lemma \ref{yinyong} (ii), one has $k_{\wz{\partial}(x_{0},x_{2})}=1$. Therefore, $\wz{\partial}(x_{0},x_{2})=(2,q-2)$. But
$P_{(1,q-1),(1,q-1)}(x_{0},x_{2})=\{x_{1}\}$ and $P_{(1,q-1),(1,q-1)}(x_{q-1},x_{1})=\emptyset$, a contradiction. Lemma \ref{yinyong} (iii) implies (\ref{2,g-2}).$\qed$

\begin{lemma}\label{CircuitQG}
For any circuit $(x_{0},x_{1},\ldots,x_{l-1})$ with $\wz{\partial}(x_{l-1},x_{0})=(1,g-1)$, there exists $i\in\{0,1,\ldots,l-2\}$ such that
$\wz{\partial}(x_{i},x_{i+1})=(1,g-1)$.
\end{lemma}
\textbf{Proof.} Suppose for the contradiction that $\wz{\partial}(x_{i},x_{i+1})=(1,q-1)$ for any $i=0,1,\ldots,l-2$. By Lemma \ref{yinyong} (iii), we have $A_{g-1,1}=A_{1,q-1}^{l-1}$. Then $A_{g-1,1}$ is a permutation matrix, a contradiction.$\qed$

\begin{lemma}\label{F(x)VGammaArcG}
$F(x)\neq V\Gamma$ if and only if every circuit of length $g$ in $\Gamma$ only contains the arcs of the same type.
\end{lemma}
\textbf{Proof.} Suppose $F(x)\neq V\Gamma$.  Assume the contrary, namely, $(x_{0},x_{1},\ldots,x_{g-1})$ is a circuit containing arcs of different types such that
$\wz{\partial}(x_{0},x_{1})=(1,g-1)$.
By (\ref{arci}) and Lemma \ref{CircuitQG},  we may assume that $\wz{\partial}(x_{1},x_{2})=(1,q-1)$ and
$\wz{\partial}(x_{g-1},x_{0})=(1,g-1)$. By the claim in Lemma \ref{Arc}, $\wz{\partial}(x_{0},x_{2})=(2,g-2)$. Since $F(x)\neq V\Gamma$ and (\ref{Arc(2)}), $\wz{\partial}(x_{g-1},x_{1})=(2,g-2)$. The fact that $x_{2}\notin F(x_{0})$ implies that $P_{(1,g-1),(1,g-1)}(x_{0},x_{2})=\emptyset$, contradicting to $x_{0}\in
P_{(1,g-1),(1,g-1)}(x_{g-1},x_{1})$.

The converse is true by Lemma \ref{F(x)VGammaArcGG}.$\qed$

\section{The proof of Theorem \ref{Main}}

In this section, we always assume that $F=\langle \Gamma_{1,g-1}\rangle$ and $x$ is a fixed vertex of $\Gamma$.

\begin{lemma}\label{gamma/f2}
If $F(x)\neq V\Gamma$, then $\Gamma/F\simeq C_{2}$.
\end{lemma}
\textbf{Proof.} Suppose for the contradiction that $\Gamma/F\simeq C_{m}$ with $m\geq3$. Choose a path $(x_{0},x_{1},x_{2},x_{3})$ such that
$\wz{\partial}(x_{0},x_{1})=\wz{\partial}(x_{1},x_{2})=(1,q-1)$ and $\wz{\partial}(x_{2},x_{3})=(1,g-1)$. Since $\partial(F(x_{0}),F(x_{2}))=2$, $k_{1,q-1}=1$ implies that $\wz{\partial}(x_{0},x_{3})=(3,l)$ for some $l$. Then there exists a shortest path $(x_{3},x_{4},y_{2},\ldots,x_{l+2},x_{0})$. By Lemma \ref{CircuitQG} and (\ref{arci}), we may assume that $\wz{\partial}(x_{3},x_{4})=(1,g-1)$. Since $\partial(F(x_{1}),F(x_{4}))=1$ and $k_{1,q-1}=1$, by (\ref{Arc(2)}), we obtain $\wz{\partial}(x_{1},x_{4})=(3,t)$ for some $t\leq l$. From $m\geq3$ and (\ref{arci}), there exists a path $(x_{4},y_{1},y_{2},\ldots,y_{t-2},x_{0},x_{1})$. Then $(x_{3},x_{4},y_{1},y_{2},\ldots,y_{t-2},x_{0})$ is a path of length $t$; and so $l\leq t$. Hence, $l=t$. By (\ref{2,g-2}), $x_{2}\in P_{(2,q-2),(1,g-1)}(x_{0},x_{3})$. Then there exists $y\in P_{(2,q-2),(1,g-1)}(x_{1},x_{4})$. From $k_{1,q-1}=1$, $\wz{\partial}(x_{2},y)=(1,q-1)$, which implies $\Gamma_{1,q-1}\in F$, a contradiction.$\qed$

\begin{prop}\label{F(x)neqVGamma}
If $F(x)\neq V\Gamma$, then $\Gamma$ is isomorphic to one of the digraphs in Theorem~\ref{Main} {\rm (i)}.
\end{prop}
\textbf{Proof.}~By Lemma \ref{gamma/f2}, $V\Gamma$ has a partition $F(x)\dot{\cup} F(x')$. Let $\Delta$ and $\Delta'$ be the subdigraphs of $\Gamma$  induced on $F(x)$ and $F(x')$, respectively. By (\ref{arci}) and $k_{1,q-1}=1$, $\sigma:F(x)\rightarrow F(x')$, $y\mapsto y'$ is an isomorphic mapping from $\Delta$ to $\Delta'$, where $y'\in \Gamma_{1,q-1}(y)$.  By Lemmas \ref{f(x)} and \ref{F(x)VGammaArcG}, $\Gamma_{r,g-r}(y)=\Delta_{r,g-r}(y)$ for each $y\in F(x)$ and $r\in\{1,2,\ldots,g-1\}$. By (\ref{arcii}), the proof of Proposition 4.3 in \cite{KSW03} implies that $\Delta$ is isomorphic to
$\Gamma_1:={\rm Cay}(\mathbb{Z}_{g}\times\mathbb{Z}_{g},\{(1,0),(0,1)\})$ or $\Gamma_2:={\rm Cay}(\mathbb{Z}_{2g},\{1,g+1\}).$ Suppose that $\tau_{i}$ is an isomorphic from $\Gamma_{i}$ to $\Delta$.

We claim that $\Delta\simeq\Gamma_{2}$. Suppose for the contrary that $\Delta\simeq\Gamma_{1}$. Write $\tau_{1}(a,b)=(a,b,0)$ and $\sigma(a,b,0)=(a,b,1)$ for each $(a,b)\in \mathbb{Z}_{g}\times\mathbb{Z}_{g}$. Let $((0,0,1),(c,d,0))$ be an arc of type $(1,q-1)$. By (\ref{2,g-2}), $\wz{\partial}_{\Gamma}((0,0,0),(c,d,0))=(2,q-2)$. Lemma \ref{f(x)} implies that $c\neq0$ and $d\neq0$. By Lemma \ref{f(x)} again, we have $(c,d,0)\in P_{(2,q-2),(g-\hat{d},\hat{d})}((0,0,0),(c,0,0))$ and $\wz{\partial}_{\Gamma}((0,0,0),(c,0,0))=\wz{\partial}_{\Gamma}((0,0,0),(0,c,0))$. By $k_{2,q-2}=1$, we have $(0,c,0)\in \Gamma_{g-\hat{d},\hat{d}}(c,d,0)$. Then $(0,c,0)\in\{(c,0,0),(c-d,d,0)\}$ by Lemma \ref{f(x)}. Hence, $c=d$.

Suppose $c=g-1$. Since $((0,0,1),(g-1,g-1,0),(0,g-1,0),(0,0,0))$ is a shortest path, $q=4$, contrary to Lemma \ref{Circuit}. Suppose $c\neq g-1$. Then $\wz{\partial}_{\Gamma}((0,0,0),(c,c+1,0))=(3,l)$ for some $l$. Pick a path $((c,c+1,0),x_{1},x_{2},\ldots,x_{l-1},(0,0,0))$. By Lemma \ref{CircuitQG} and ({\ref{arci}}), we may assume that $\wz{\partial}_{\Gamma}((c,c+1,0),x_{1})=(1,g-1)$. By (\ref{Arc(2)}), we have $\wz{\partial}_{\Gamma}((0,0,1),x_{1})=(3,t)$ for some $t\leq l$. Since $F(x)\neq V\Gamma$, $k_{1,q-1}=1$ implies that there exists a path $(x_{1},y_{1},y_{2},\ldots,y_{t-2},(0,0,0),(0,0,1))$. Then $((c,c+1,0),x_{1},y_{1},y_{2},\ldots,y_{t-2},(0,0,0))$ is a path of length $t$; and so $l\leq t$. Hence $l=t$. By (\ref{2,g-2}) and $x_{1}\in V\Delta$, one has $(c,c,0)\in P_{(2,q-2),(1,g-1)}((0,0,0),(c,c+1,0))$ and $P_{(2,q-2),(1,g-1)}((0,0,1),x_{1})=\emptyset$ in $\Gamma$, a contradiction. Therefore, our claim is valid.

Write $\tau_{2}(a)=(a,0)$ and $\sigma(a,0)=(a,1)$ for each $a\in \mathbb{Z}_{2g}$. Let $((a,1),(a+k_{a},0))$ be an arc of type $(1,q-1)$. Then $k_{a}\neq0$. By (\ref{2,g-2}), $\wz{\partial}_{\Gamma}((a,0),(a+k_{a},0))=(2,q-2)$. By Lemma \ref{f(x)}, $\wz{\partial}_{\Delta}((a,0),(a+k_{a},0))\neq (t,g-t)$ for any $t\in\{1,2,\dots,g-1\}$. Since $\bigcup_{1\leq t\leq g-1}\Delta_{t,g-t}(a,0)=V\Delta\setminus\{(a,0),(a+g,0)\}$, one has $k_{a}=g$. Then,
$\Gamma\simeq\textrm{Cay}(\mathbb{Z}_{4}\times\mathbb{Z}_{g},\{(0,1),(1,0),(2,1)\})$ and the result holds by Proposition \ref{G_Wdrdg}.$\qed$

\begin{lemma}\label{F(x)VGammaArc}
If $F(x)=V\Gamma$, then $p_{(1,g-1),(1,g-1)}^{(1,q-1)}=2$.
\end{lemma}
\textbf{Proof.}~By Lemma \ref{F(x)VGammaArcGG}, there exists a circuit of length $g$ with different types of arcs. Let $C:=(x_{0},x_{1},\ldots,x_{g-1})$ be such a circuit with the minimum number of arcs of type $(1,g-1)$. Suppose $C$ contains $t$ arcs of types $(1,g-1)$. Lemma \ref{CircuitQG} implies that $t\geq2$. By (\ref{arci}), we may assume that $\wz{\partial}(x_{i},x_{i+1})=(1,g-1)$ for $0\leq i\leq t$. We claim
that $\wz{\partial}(x_{0},x_{2})=(1,q-1)$. Suppose not. By the claim in Lemma \ref{Arc} and (\ref{Arc(2)}), we have $\wz{\partial}(x_{g-1},x_{1})=\wz{\partial}(x_{0},x_{2})=(2,g-2)$. Since $x_{0}\in P_{(1,q-1),(1,g-1)}(x_{g-1},x_{1})$, there exists $x_{1}'\in P_{(1,q-1),(1,g-1)}(x_{0},x_{2})$. The circuit
$C':=(x_{0},x_{1}',x_{2},\ldots,x_{g-1})$ contains just $t-1$ arcs of type $(1,g-1)$, a contradiction. Thus, our claim is valid. It follows that $p_{(1,q-1),(g-1,1)}^{(1,g-1)}=1$. By Lemma \ref{yinyong} (i), the desired result holds.$\qed$

Let $H=\langle\Gamma_{1,q-1}\rangle$ and $H(x_{0,0}),H(x_{0,1}),\ldots,H(x_{0,s-1})$ be all pairwise distinct vertices of $\Gamma/H$. Since $q<g$, the subdigraph induced on each $H(x_{0,j})$ is a circuit of length $q$ with arcs of type $(1,q-1)$, say $(x_{0,j},x_{1,j},\ldots,x_{q-1,j})$. It follows that $s\geq2$.

\begin{prop}\label{F(x)Gamma}
If $F(x)=V\Gamma$, then $\Gamma$ is isomorphic to one of the digraphs in Theorem~\ref{Main} {\rm (ii)}.
\end{prop}
\textbf{Proof.}~Suppose $\partial(H(x_{0,0}),H(x_{0,1}))=1$. By (\ref{arci}), we may assume that $\wz{\partial}(x_{0,0},x_{0,1})=(1,g-1)$. By Lemma \ref{F(x)VGammaArc}, one has $\wz{\partial}(x_{0,1},x_{1,0})=(1,g-1)$, which implies that $\partial(H(x_{0,1}),H(x_{0,0}))=1$. Since $F(x)=V\Gamma$, $\Gamma/H$ is a connected undirected graph. By $k_{1,g-1}=2$, $\Gamma/H$ is an undirected circuit of length $s$. Suppose $s=2$. Pick $y\in \Gamma_{1,g-1}(x_{0,1})\setminus\{x_{1,0}\}$. Then $y=x_{i,0}$ for some $i\geq2$, and $(x_{0,1},y,x_{i+1,0},\ldots,x_{q-1,0},x_{0,0})$ is a path of length $q-i+1$ from $x_{0,1}$ to $x_{0,0}$,   contrary to the fact $\partial(x_{0,1}, x_{0,0})=g-1$. Hence, $s\geq3$.

Let $(H(x_{0,0}),H(x_{0,1}),\ldots,H(x_{0,s-1}))$ be an undirected circuit. By (\ref{arci}), we may assume that $(x_{0,0},x_{0,1},\ldots,x_{0,s-1})$ is a path with arcs of type $(1,g-1)$. By Lemma \ref{F(x)VGammaArc},
$(x_{0,j},x_{0,j+1},x_{1,j},x_{1,j+1},x_{2,j},\ldots,x_{q-1,j},x_{q-1,j+1})
$
is a circuit with arcs of type $(1,g-1)$ for any $j=0,1,\ldots,s-2$. Therefore, there exists $k\in \{1,2,\ldots, q\}$  such that $\wz{\partial}(x_{0,s-1},x_{q-k+1,0})=(1,g-1)$, where the first subscription of $x$ are taken modulo $q$. By Lemma \ref{F(x)VGammaArc} again, $\wz{\partial}(x_{i,s-1},x_{i-k+1,0})=\wz{\partial}(x_{i-k+1,0},x_{i+1,s-1})=(1,g-1)$ for each $i$. Since
$(x_{0,0},x_{0,1},\ldots,x_{0,s-1},x_{q-k+1,0},x_{q-k+2,0},\ldots,x_{q-1,0})
$
is a circuit of length $s+k-1$ with different types of arcs,  By Lemma \ref{Circuit} we get $s+k-1>q$. By Proposition \ref{QSK_Wdrdgiff}, the desired result follows.$\qed$

Combining Propositions \ref{F(x)neqVGamma} and \ref{F(x)Gamma}, we complete the proof of Theorem \ref{Main}.

\section*{Acknowledgement}
This research is supported by NSFC(11271047, 11301270, 11371204) and the Fundamental Research Funds for the Central University of China.

\end{CJK*}

\end{document}